\documentclass[letterpaper, 10 pt, journal, twoside]{IEEEtran}

\IEEEoverridecommandlockouts        

\usepackage{cite}
\usepackage{amsmath,amssymb,amsfonts}
\usepackage{algorithmic}
\usepackage{hyperref}
\usepackage{textcomp}
\usepackage{graphicx,color}
\usepackage{mathrsfs}
\usepackage[vlined,ruled]{algorithm2e}
\usepackage{subfigure}
\usepackage{url}
\usepackage{color}
\usepackage{dsfont}
\usepackage{bbm}
\usepackage{booktabs}
\usepackage{array}
\usepackage[table]{xcolor} 
\usepackage{yfonts}

\newtheorem{theorem}{Theorem}[section]

\newtheorem{corollary}[theorem]{Corollary}

\newtheorem{remark}{Remark}

\newtheorem{example}{Example}

\newcommand{\map}[3]{#1: #2 \rightarrow #3}

\newcommand{\subscr}[2]{{#1}_{\textup{#2}}}

\newcommand{\until}[1]{\{1,\dots,#1\}}

\newcommand{\Ker}{\operatorname{Ker}}

\newcommand{\Image}{\operatorname{Im}}

\newcommand{\real}{\mathbb{R}}

\newcommand{\transpose}{\mathsf{T}} 
\newcommand{\E}{\operatorname{\mathbb{E}}}

\newcommand{\1}{\mathds{1} }

\DeclareSymbolFont{bbold}{U}{bbold}{m}{n}
\DeclareSymbolFontAlphabet{\mathbbold}{bbold}

\newcommand{\vect}[1]{\mathbbold{#1}}
\newcommand{\vectorones}[1][]{\vect{1}_{#1}}

\newcommand{\R}{\mathbb{R}}

\newcommand{\N}{\mathbb{N}}

\newcommand\oprocendsymbol{\hbox{$\square$}}
\newcommand\oprocend{\relax\ifmmode\else\unskip\hfill\fi\oprocendsymbol}

\renewcommand{\theenumi}{(\roman{enumi})}

\newcommand*{\QEDA}{\hfill\ensuremath{\blacksquare}}%

\graphicspath{{img/}}

\begin{document}
\title{Data-Driven Minimum-Energy Controls\\ for Linear Systems} \author{Giacomo Baggio, Vaibhav Katewa, and Fabio
  Pasqualetti \thanks{This material is based upon work supported in
    part by ARO 71603NSYIP. Giacomo Baggio, Vaibhav Katewa and Fabio
    Pasqualetti are with the Department of Mechanical Engineering,
    University of California at Riverside,
    \{\href{mailto:gbaggio.ucr.edu}{\texttt{gbaggio}},
    \href{mailto:vkatewa@engr.ucr.edu}{\texttt{vkatewa}},
    \href{mailto:fabiopas@engr.ucr.edu}{\texttt{fabiopas\}@engr.ucr.edu.}}}}
\maketitle

\thispagestyle{empty}

\begin{abstract}
  In this paper we study the problem of computing minimum-energy
  controls for linear systems from experimental data. The design of
  open-loop minimum-energy control inputs to steer a linear system
  between two different states in finite time is a classic problem in
  control theory, whose solution can be computed in closed form using
  the system matrices and its controllability Gramian. Yet, the
  computation of these inputs is known to be ill-conditioned,
  especially when the system is large, the control horizon long, and
  the system model uncertain. Due to these limitations, open-loop
  minimum-energy controls and the associated state trajectories have
  remained primarily of theoretical value. Surprisingly, in this paper
  we show that open-loop minimum-energy controls can be learned
  exactly from experimental data, with a finite number of control
  experiments over the same time horizon, without knowledge or
  estimation of the system model, and with an algorithm that is
  significantly more reliable than the direct model-based
  computation. These findings promote a new philosophy of controlling
  large, uncertain, linear systems where data is abundantly available.
\end{abstract}

\section{Introduction}\label{sec: introduction}
\IEEEPARstart{C}{onsider}  the discrete-time linear time-invariant system
\begin{align}\label{eq: system}
  x(t+1) = A x(t) + Bu(t),
\end{align}
where, respectively, $A\in\R^{n\times n}$ and $B\in\R^{n\times m}$
denote the system and input matrices, and $\map{x}{\N}{\R^n}$ and
$\map{u}{\N}{\R^m}$ describe the state and input of the system. For a
control horizon $T \in \N$ and a desired state $\subscr{x}{f}$, the
minimum-energy control problem asks for the input sequence
\mbox{$u(0), \dots, u(T-1)$} with minimum energy that steers the state from
$x_{0}$ to $\subscr{x}{f}$ in $T$ steps, and it can be formulated as
\begin{align}\label{eq: min energy problem 1}
  \begin{array}{ll}
    \min\limits_{u} & \sum\limits_{t=0}^{T-1} \| u(t)\|_2^2 ,\\[1em]
    \,\text{s.t.} 
         & x(t+1) = A x(t) + B u(t), \\[.5em]
         & x(0) =x_{0}, \ x(T) = \subscr{x}{f} .
  \end{array}
\end{align}
As a classic result \cite{TK:80}, the minimization problem \eqref{eq:
  min energy problem 1} is feasible if and only if
$(\subscr{x}{f} - A^T x_0) \in \Image (W_T)$, where
\begin{align}\label{eq: Gramian}
  W_T = \sum_{t = 0}^{T-1} A^t B B^\transpose (A^\transpose)^t
\end{align}
is the $T$-steps controllability Gramian and $\Image (W_T)$ denotes
the image of the matrix $W_T$. Further, the solution to \eqref{eq: min
  energy problem 1} is
\begin{align}\label{eq: min energy input}
  u^* (t) = B^\transpose (A^\transpose)^{T-t-1} W_T^{\dag}
  (\subscr{x}{f} - A^T x_0) ,
\end{align}
where $W_T^\dag$ is the Moore--Penrose pseudoinverse of $W_T$
\cite{AB-TNEG:03}.

The controllability Gramian \eqref{eq: Gramian} and the minimum-energy
control input \eqref{eq: min energy input} identify fundamental
control limitations for the system \eqref{eq: system}, and have been
extensively used to solve design \cite{SZ-FP:16a}, sensor and actuator
placement \cite{THS-FLC-JL:16}, and control 
problems \cite{FP-SZ-FB:13q} for systems and
networks. However, besides their theoretical value, the optimal
control input \eqref{eq: min energy input} is rarely used in practice
or even computed numerically because (i) it relies on the perfect
knowledge of the system dynamics, (ii) its performance is not robust
to model uncertainties, and (iii) the controllability Gramian is
typically ill-conditioned, especially when the system is large
\cite{FP-SZ-FB:13q,DCS-YZ:02b}. This implies that the control sequence
\eqref{eq: min energy input} is numerically difficult to compute, and
that its implementation leads to errors \cite{JS-AEM:13}. To the best
of our knowledge, efficient and numerically reliable methods to
compute minimum-energy control inputs are still lacking.

\noindent
\textbf{Paper contributions.} This paper features two main
  contributions. First, we show that minimum-energy control inputs for
  linear systems can be computed from data obtained from control
  experiments with non-minimum-energy inputs, and without knowledge or
  estimation of the system matrices. Thus, optimal inputs can be
  learned from \emph{non-optimal} ones, and we provide three different
  expressions for doing so. Surprisingly, we also establish that a
  \emph{finite} number of non-optimal control experiments is always
  sufficient to compute minimum-energy control inputs towards any
  reachable state. Second, we show that the data-driven computation of
  minimum-energy inputs is numerically as reliable as the computation
  of the inputs based on the exact knowledge of the system matrices,
  and substantially more reliable than using the closed-form
  expression based on the Gramian. Further, as minor contributions, we
  (i) derive bounds on the number of required control experiments as a
  function of the dimension of the system, number of control inputs,
  and length of the control horizon, (ii) discuss the effect of noisy
  data on the data-driven expressions, and (iii) extend our
  data-driven framework to the case of output measurements. 

Our results suggest the tantalizing hypothesis that several optimal
control problems can be solved efficiently and reliably using a
combination of data-driven algorithms and system properties (in our
setup, linearity of the dynamics), even when the system model is
uncertain or unknown.

\noindent
\textbf{Related work.}  Several works investigate the problem of
  estimating optimal controls for linear systems from input-output
data. The classic model-based approach \cite{MG:05} consists of (i)
identifying a model of the system from the available data, and (ii)
using the estimated model to design the optimal control
inputs. Data-driven algorithms have been proposed in
\cite{KF-MW:95,GS-RES:00,WA-DK-BJ-RM-MS:05,GRGS-ASB-CL-LC:19}
for the LQR/LQG problem. In particular, the approach pursued in these
papers relies on the estimation of the Markov parameters of the
system, thereby bypassing the identification step of the model-based
approach. Differently from the above approaches, in this paper we
focus on computing open-loop minimum-energy inputs from
experimental data, without reconstructing the system matrices and
where the experiments use arbitrary control inputs. To the best of our
knowledge, this paper addresses a novel problem and provides new and
numerically more reliable expressions for the computation of
minimum-energy control inputs.

\section{Learning minimum-energy control inputs}\label{sec:
  learning}
  
In vector form, the minimum-energy control problem asks to find the
minimum-norm solution to the following equation:
\begin{align*}
  \subscr{x}{f} =
  A^T x_0 + 
  \underbrace{
  \begin{bmatrix}
    B & AB & \cdots & A^{T-1}B
  \end{bmatrix}}_{C_T}
                      u,
\end{align*}
where the vector $u\in\R^{mT}$ contains the control inputs over the control
horizon $[0,T-1]$, namely $u=[u(T-1)^{\transpose}\, \cdots\,\, u(0)^{\transpose}]^{\transpose}$, and $C_T$ denotes the $T$-steps controllability
matrix.\footnote{To simplify the technical treatment and without
  compromising generality, we assume that $\subscr{x}{f}$
  is reachable in $T$-steps,~i.e., $\!(\subscr{x}{f}- A^{T}x_0 )\in \Image(C_T)$.} Then, if the controllability
matrix $C_T$ is known, the minimum-energy control input to reach
$\subscr{x}{f}$ is
\begin{align}\label{eq: model-based}
  u^* = C_T^\dag (\subscr{x}{f} - A^T x_{0}).
\end{align}
Instead of using \eqref{eq: model-based}, in this paper we aim to
compute minimum-energy control inputs leveraging a set of $N$ control
experiments and assuming that the system matrices, and thus the
controllability matrix, are not available. 
The $i$-th control experiment
  consists of applying the input sequence $u_i$ to
  \eqref{eq: system}, and measuring the system state at
  time $T$, namely $x_{i}$, where
\begin{align}\label{eq: data relation}
  x_i = A^Tx_{0} + C_T u_i. 
\end{align}
We remark that the inputs $u_i$ are arbitrary and not necessarily of
minimum-norm. In vector form, the available data is
\begin{align}\label{eq: X and U}
  X = 
  \begin{bmatrix}
    x_1 &\! \cdots\! & x_N
  \end{bmatrix},\ \text{ and }\ \ 
                   U  = 
                   \begin{bmatrix}
                     u_1 &\! \cdots  \!  & u_N 
                   \end{bmatrix},          
\end{align}
where $x_{i}$ is the state at time $T$  with input $u_{i}$ as in \eqref{eq: data relation}.\footnote{While the full state trajectory could be measured \cite{CDP-PT:19}, here we show that measuring the final state is sufficient to compute minimum-energy~inputs.}

\subsection{Data-driven minimum-energy controls}

Because we only rely on the experimental data $(X,U)$ to learn the
minimum-energy control input to reach a desired state, we postulate
that such input can be computed as a linear combination of the inputs
$U$. Thus, we formulate and study the following constrained
minimization problem:
\begin{align}\label{eq: min energy problem}
  \begin{array}{lcl}
    \alpha^* = &\arg\min\limits_{\alpha} &  \| U \alpha \|_2^2 ,\\[1em]
    & \text{s.t.} 
                         & \subscr{x}{f} = X \alpha,
  \end{array}
\end{align}
where $\alpha\in\real^{N}$ is the optimization variable.
As we show in Theorem \ref{thm: feasibility}, a first data-driven
expression for the minimum-energy control input derives from a
solution to \eqref{eq: min energy problem}. We start with the
expression of the minimum-energy control input for the case $x_0 =0$,
and we postpone the general case $x_0 \neq 0$ to Remark \ref{remark:
  nonzero initial state}. Let $\Image(M)$ and $\Ker(M)$ denote the
range-space and the null-space of the matrix $M$, respectively. With a
slight abuse of notation, we write $K~=~\Image(A)$ (resp.
$K = \Ker(A)$) to say that $K$ is a basis of $\Image(A)$
(resp. $\Ker(A)$). A matrix is full row rank if the dimension of its
range-space equals the number of its rows.

\begin{theorem}{\bf \emph{(Data-driven minimum-energy
      control inputs when $x_0 = 0$)}}\label{thm: feasibility}
  If the matrix $U$ in \eqref{eq: X and U} is full row rank, then, for
  any final state $\subscr{x}{f}$, the minimum-energy input~equals
  \begin{align}\label{eq: U min 1}
    u^* = (I - UK (UK)^\dag) U X^\dag \subscr{x}{f},
  \end{align}
  where $K = \Ker(X)$ and $X$ is as in \eqref{eq: X and U}.
\end{theorem}
\begin{IEEEproof} We first show that \eqref{eq: min energy problem}
    is feasible, and that $u^{*}=U\alpha^{*}$. Notice that, because
  $U$ is full row rank, there exists $\alpha^*$ such that
  $u^* = U \alpha^*$, where $u^*$ is the minimum-energy control input
  to reach $\subscr{x}{f}$. Additionally, $\alpha^*$ satisfies the
  constraint in \eqref{eq: min energy problem} because
  $X \alpha^* = C_T U \alpha^* = C_T u^* =
  \subscr{x}{f}$. Finally, because $u^*$ is unique \cite{TK:80},
  $\alpha^*$ is also a solution to \eqref{eq: min energy problem}, and
  its computation is equivalent to computing the 
  input~$u^*$. 
  
  To compute $\alpha^*$ we solve the constraint $\subscr{x}{f} = X \alpha$ and substitute it in
  the cost function. Namely, $\alpha = X^\dag \subscr{x}{f} - K w$,
  where $K = \Ker(X)$ and $w$ is an arbitrary vector. Equating to zero
  the derivative of the cost function with respect to $w$, we obtain
  $
   w^{*} = (UK)^\dag U X^\dag
    \subscr{x}{f}.
  $
This implies that $\alpha^* = X^\dag \subscr{x}{f} - K w^*$, from which \eqref{eq: U min 1} follows by letting~$u^* = U \alpha^*$.~
\end{IEEEproof}

Theorem \ref{thm: feasibility} provides an expression of the
minimum-energy control input, which only uses data originated from a
set of control experiments, and does not require the knowledge of the
system matrices. Importantly, Theorem \ref{thm: feasibility} shows
that minimum-energy control inputs can be directly computed based on a
number of control experiments with arbitrary, thus not minimum-energy,
inputs. Further, Theorem \ref{thm: feasibility} assumes that $U$ is
full row rank, which guarantees the computation of the minimum-energy
input for any final state $\subscr{x}{f}$. When $U$ is not full row
rank but $u^* \in \Image (U)$, the minimum-energy control input can
still be computed as in Theorem \ref{thm: feasibility}. Instead, when
$u^* \not\in \Image (U)$, the minimum-energy input cannot be computed
as a (linear) combination of the experimental data \eqref{eq: X and
  U}. In this case, the data-driven input \eqref{eq: U min 1}
reaches the desired final state $\subscr{x}{f}$, if
$\subscr{x}{f}\in \Image(X)$, or the final state
$\subscr{\tilde{x}}{f} \in \Image(X)$ that is closest to
$\subscr{x}{f}$, if $\subscr{x}{f}\not\in \Image(X)$.  To see this,
let $u^*$ be as in \eqref{eq: U min 1} and note~that
\begin{align*}
  \subscr{\tilde{x}}{f} &= C_{T} u^* = C_{T} (I - UK
                          (UK)^\dag) U  X^\dag \subscr{x}{f} \\
                        &= C_{T} U X^\dag \subscr{x}{f} - \underbrace{C_{T} UK
                          (UK)^\dag U  X^\dag \subscr{x}{f}}_{=0
                          \text{ because }
                          C_{T} UK = XK = 0} = X X^\dag \subscr{x}{f},
\end{align*}
which shows that $\subscr{\tilde{x}}{f}$ is the orthogonal projection
of $\subscr{x}{f}$ onto $\Image(X)$.  This in particular implies that
the error $\|\subscr{x}{f}-\subscr{\tilde{x}}{f}\|_2$ is
non-increasing in the number of experiments $N$, and it vanishes when
the experimental data satisfies $\subscr{x}{f} \in \Image(X)$.
Finally, Theorem \ref{thm: feasibility} can also be used to quantify
the number of experiments needed to compute~minimum-energy~inputs.

\begin{corollary}{\bf \emph{(Required number of control experiments to
      compute minimum-energy inputs)}}\label{corollary: number
    experiments}
  Let $n$ be the dimension of the system, $m$ the number of inputs,
  $T$ the control horizon, and $N$ the number of control
  experiments. Then,
  \begin{enumerate}
  \item $N \ge n$ is necessary to compute minimum-energy control
    inputs towards any arbitrary final state $\subscr{x}{f}$;

  \item $N =  mT$ is sufficient to compute minimum-energy control
    inputs towards any arbitrary final state $\subscr{x}{f}$, provided
    that the inputs $u_i$ are linearly independent.
  \end{enumerate}
\end{corollary}
\begin{IEEEproof}
  \textit{(Necessity)} Assume by contradiction that the number of
  experiments is strictly less than $n$. Then, $\text{Rank}(X) < n$,
  and there exists $\subscr{x}{f} \not\in \Image(X)$. Then, the
  minimization problem \eqref{eq: min energy problem} is infeasible,
  and the minimum-energy control input cannot be computed from the
  inputs~$U$.
  
  \textit{(Sufficiency)} Let the experimental inputs be linearly
  independent. Then, $U$ is invertible and, for any $\subscr{x}{f}$,
  there exists a solution $\alpha^*$ such that $u^* = U\alpha^*$. This
  shows that the minimum-energy input can be computed from the data.
\end{IEEEproof}

Corollary \ref{corollary: number experiments} characterizes the number
of control experiments that are required to compute minimum-energy
control inputs from experimental data. In particular, as few as $n$
experiments are needed, in which case the experiments must contain $n$
linearly independent minimum-energy control inputs, and as many as
$mT$ experiments are sufficient, in which case the control inputs can
be selected arbitrarily provided that they form a linearly independent
set of vectors. This also shows that optimal control inputs can be
learned from a \emph{finite} number of non-optimal control
inputs. 

\begin{figure}
  \includegraphics[scale=1]{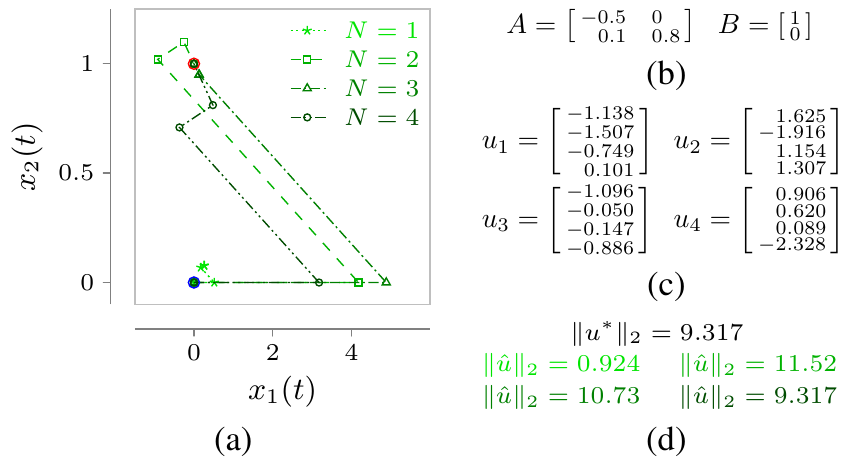}
  \caption{Fig.~\ref{Fig:2d-trajectories}(a) illustrates the state
    trajectory $x(t)=[x_{1}(t) \ x_{2}(t)]^{\transpose}$, $t\in[0,T]$,
    of the system described by matrices $A$, $B$ as in
    Fig.~\ref{Fig:2d-trajectories}(b) and driven by the data-driven
    input $\hat{u}$ \eqref{eq: U min 1}, for four values of the number
    of experiments $N$. We choose $T=4$ and denote the initial state
    $x_0=[0 \ 0]^{\transpose}$ and final state
    $\subscr{x}{f}=[0 \ 1]^{\transpose}$ with a blue and red circle,
    respectively.  The other markers correspond to the values of the
    trajectories in the interval $[0,T]$. The data-driven input
    $\hat{u}$ has been computed using input data as in $u_{i}$ as
    Fig.~\ref{Fig:2d-trajectories}(c), and $x_{i}=C_{T}u_{i}$.
    Fig.~\ref{Fig:2d-trajectories}(d) shows the norm of the
    minimum-energy input and the norm of the data-driven input
    $\hat{u}$ \eqref{eq: U min 1} as $N$ varies~(color-coded).}
  \label{Fig:2d-trajectories}
\end{figure}

  \begin{example}{\bf \emph{(Data-driven control inputs when
        $N\le mT$)}}\label{ex:2d-sys}
    We consider a two-dimensional system with matrices $A$ and $B$ as in
    Fig.~\ref{Fig:2d-trajectories}(b), control horizon $T=4$, initial state $x_0=[0 \ 0]^{\transpose}$, and final state $\subscr{x}{f}=[0 \ 1]^{\transpose}$. 
    We vary the
      number of control experiments $N$ from $1$ to $4$, where the
      inputs are as in Fig.~\ref{Fig:2d-trajectories}(c). For each
      number of control experiments, we compute the data-driven input
      \eqref{eq: U min 1}, and report the corresponding state
      trajectory and norm in Fig.~\ref{Fig:2d-trajectories}(a) and
      Fig.~\ref{Fig:2d-trajectories}(d), respectively. Notice that,
      when $N = 1$, the data-driven input does not steer the system
      state to $\subscr{x}{f}$. Instead, for $N=2,3,4$ the state
      trajectory reaches~$\subscr{x}{f}$. Finally, the data-driven
      input has minimum norm only when $N=4$.~\oprocend
\end{example}

\begin{remark}{\bf \emph{(Geometric properties of \eqref{eq: U min
        1})}}\label{remark: geometric properties}
  Several geometric properties of \eqref{eq: U min 1} can be
  highlighted. First, $UK = \Ker (C_T)$ when $U$ is full row rank. In
  fact, $C_T UK = X K = 0$, showing that
  $\Image(UK) \subseteq \Ker(C_T)$. Further, if $C_T u = 0$ and
  $u = U \alpha$, then, $X \alpha = C_T U \alpha = C_T u = 0$, showing
  that $\alpha \in \Image(K)$ and $\Ker(C_T)\subseteq
  \Image(UK)$. Thus, $\Image(UK) = \Ker( C_T)$ when $U$ is full row
  rank. Second, $I - UK (UK)^\dag$ is the orthogonal projection onto
  the kernel of $(UK)^\transpose$ and, consequently,
  $u^* = (I - UK (UK)^\dag) U X^\dag \subscr{x}{f}$ is orthogonal to
  $\Ker(C_T)$. This~is~expected, because $u^*$ is the minimum-energy
  control input to reach the state $\subscr{x}{f}$. \oprocend
\end{remark}

\subsection{An alternative expression of minimum-energy controls}

In this subsection, we present a different optimization problem that can
be used to derive an equivalent expression of the data-driven
minimum-energy control input \eqref{eq: U min 1}. Specifically, we consider the
following problem, which encodes the problem of
estimating the controllability matrix from data:
\begin{align}\label{eq:ctrb_mat_estimate}
  \begin{array}{ll}
    C_{T}^* = \arg\min\limits_{C} &  \| X - CU \|_F^2,
  \end{array}
\end{align}
where $\|\cdot\|_F$ denotes the Frobenius norm of a matrix. The above
problem has a unique solution, which equals
\mbox{$C_{T}^{*} = XU^{\dag}$}. Notice that the minimization problem
\eqref{eq:ctrb_mat_estimate} returns an estimate of the
controllability matrix, which can be used to compute the input as
$\hat{u}=(C^{*}_{T})^{\dag}\subscr{x}{f}=(XU^{\dag})^{\dag}\subscr{x}{f}$. We
next show that $\hat{u}$ coincides with the control input \eqref{eq: U
  min 1}.

\begin{theorem}{\bf \emph{(Equivalent expressions of data-driven
      minimum-energy inputs)}}\label{thm: equivalent expression}
  Let $X$ and $U$ be as in \eqref{eq: X and U}. Then,
  \begin{align}\label{eq: alternative expression}
    (I - UK (UK)^\dag) U X^\dag  \subscr{x}{f} = (XU^\dag)^\dag  \subscr{x}{f}.
  \end{align}
\end{theorem}\medskip
\begin{IEEEproof}
  We show that $(XU^{\dag})^{\dag} = (I - UK (UK)^\dag) U
  X^\dag$. That is, we show that
    $(I-UK(UK)^{\dagger})UX^{\dagger}$ satisfies the four conditions
    \cite{AB-TNEG:03} defining the Moore--Penrose pseudoinverse of
    $XU^{\dagger}$. To this aim, let $K = I-X^{\dag} X$. Since
    $P = I - UK (UK)^\dag$ is the orthogonal projection onto
  $\text{Ker}((UK)^{\transpose})$,
  \begin{align}\label{eq:eq_pf_prop_1}
    (UK)^{\transpose} P = 0
    \overset{P=P^{\transpose}}{\Longrightarrow} PUK = 0 \Rightarrow PU
    X^{\dag} X = PU.
  \end{align}
  Because $X = C_TU$, we have
  $\text{Ker}(U) \subseteq \text{Ker}(X)$. Since $I-U^{\dag} U$ is the
  orthogonal projection onto $\text{Ker}(U)$, we have
  \begin{align}\label{eq:eq_pf_prop_2}
    X (I-U^{\dag} U) = 0 \Rightarrow X U^{\dag} U = X.
  \end{align}
  Further, using $XK = 0$, we obtain
  \begin{align}\label{eq:eq_pf_prop_3}
    X U^{\dag} (I-P) \! = \!X U^{\dag} UK (UK)^{\dag}
    \overset{\eqref{eq:eq_pf_prop_2}}{=} XK (UK)^{\dag}\! =\! 0.
  \end{align}
  Finally, since $I-U U^{\dag}$ denotes the orthogonal projection
  onto $\text{Ker}(U^{\transpose})$, and $UK (UK)^\dag$ the
  orthogonal projection onto
  $\text{Im}(UK) \subseteq \text{Im}(U) \perp
  \text{Ker}(U^{\transpose})$, we have
  \begin{align} \label{eq:eq_pf_prop_4}
    &\ \ \ \ \ \, UK (UK)^\dag = I-P=0 \nonumber \\
    &\Rightarrow\  (I-P)(I-U U^{\dag}) = [(I-P)(I-U U^{\dag})]^{\transpose} \nonumber \\
    & \Rightarrow\ U U^{\dag} P = P U U^{\dag},
  \end{align}
  where the last implication follows because $I-P$ and
  $I-U U^{\dag}$ are symmetric. To conclude, we show that
  $PUX^{\dag} = (XU^{\dag})^\dag$ by proving the four
    Moore--Penrose conditions \cite{AB-TNEG:03}: 
  \begin{enumerate} 
  \item \phantom{I} $PUX^{\dag} X U^{\dag} PUX^{\dag} 
    \overset{\eqref{eq:eq_pf_prop_1}}{=}  P U U^{\dag} PUX^{\dag}
    \overset{\eqref{eq:eq_pf_prop_4}}{=} P^2 U  U^{\dag}\cdot U X^{\dag}\- =
    PUX^{\dag}$; 

  \item
    $XU^{\dag} PUX^{\dag} X U^{\dag}
    \overset{\eqref{eq:eq_pf_prop_1}}{=} X U^{\dag} PU U^{\dag} = X
    U^{\dag} U U^{\dag} - X U^{\dag} (I-P)U U^{\dag}
    \overset{\eqref{eq:eq_pf_prop_3}}{=} X U^{\dag}$;

  \item
    $X U^{\dag} PUX^{\dag} = X U^{\dag} UX^{\dag} - X U^{\dag}(I-P)UX^{\dag} \hspace*{2pt}\overset{\eqref{eq:eq_pf_prop_2},\, 
      \eqref{eq:eq_pf_prop_3}}{=} XX^{\dag} = (XX^{\dag})^{\transpose}$;

  \item
    $PUX^{\dag} X U^{\dag} \overset{\eqref{eq:eq_pf_prop_1}}{=} PU
    U^{\dag} \overset{\eqref{eq:eq_pf_prop_4}}{=} U U^{\dag} P = (PU
    U^{\dag})^{\transpose}$.
  \end{enumerate}
\end{IEEEproof}

\subsection{An asymptotic expression of minimum-energy controls}

The minimization problem \eqref{eq:ctrb_mat_estimate} reconstructs the
forward controllability matrix $C_T$, from which minimum-energy
control inputs can be derived by subsequently computing $C_T^\dag$. To
avoid the computation of $C_T^\dag$ and obtain a potentially simpler
expression, we next consider the problem of directly estimating
$C_T^\dag$ from the experimental data:
\begin{align}\label{eq: inverse dynamics}
  \begin{array}{ll}
    M^* = \arg\min\limits_{M} &  \| MX - U \|_F^2 .
  \end{array}
\end{align}
Notice that the latter problem is equivalent to estimating the inverse
map from $X$ to $U$, and it is typically more difficult than the
problem of estimating the map from $U$ to $X$. In fact, while the
forward map is unique, the inverse map is typically not.\footnote{In particular, the inverse map is not unique whenever $mT>n$.} Further, the control
input $M^* \subscr{x}{f}$ obtained by solving the minimization problem
\eqref{eq: inverse dynamics} is not guaranteed to be of minimum norm
and to steer the system to $\subscr{x}{f}$, as these constraints do
not appear in the minimization problem. In what follows, we say that a sequence of random matrices $\{X_{n}\}_{n\in\N}$ converges almost surely (a.s.) to a matrix $X$, and denote it with $X_{n}\xrightarrow{\text{a.s.}}X$, if $\mathrm{Pr}(\lim_{n\to \infty} X_{n}=X)=1$.
\begin{theorem}{\bf \emph{(Asymptotically equivalent expression to
      \eqref{eq: U min 1})}}\label{thm: equivalent expression 2}
  Let $X$ and $U$ be as in \eqref{eq: X and U}. The unique solution to the
  minimization problem \eqref{eq: inverse dynamics} is
  \begin{align}\label{eq: solution M}
     M^* = U X^\dag ,
  \end{align}
  and the corresponding control input can be written as
  \begin{align}\label{eq: solution M input}
    \hat{u} =M^*\subscr{x}{f} = U X^\dag \subscr{x}{f}.
  \end{align}
  Further, if $X$ is full row rank, then
  $C_T M^* \subscr{x}{f} = \subscr{x}{f}$. That is, the control
  $\hat{u}$ steers the system from $x_0 = 0$ to
  $x(T)~=~\subscr{x}{f}$.
  Finally, if
  the entries of $U$ are i.i.d. random variables with zero
  mean and nonzero finite variance, then $U X^\dag  \xrightarrow{\text{a.s.}} C_T^\dag$ as $N\to\infty$.
  That is, as the number of control experiments increases, the input
  $\hat{u}$ converges a.s. to the optimal input~$u^*$.
\end{theorem}
\begin{IEEEproof}
  The expression \eqref{eq: solution M} follows from the properties of
  the Moore--Penrose pseudoinverse. For the second claim,
 we note that $
    C_T \hat{u} = C_T U X^{\dag}
    \subscr{x}{f} = XX^{\dag} \subscr{x}{f}
    = \subscr{x}{f}$,
  where we have used that $X$ is full row rank and $X = C_T U$. To
  prove the third statement, let $N\to \infty$, and let the control
  experiments be chosen so that the entries of $U$ are i.i.d. random
  variables with zero mean and finite variance $\sigma^2$. Let
  $U_{ij}$ denote the $(i,j)$-th entry of $U$, and observe that the
  $(i,j)$-th entry of $\frac{1}{N}UU^\transpose$ equals
  $\frac{1}{N}\sum_{k=1}^{N} U_{ik}U_{jk}$. Because
  $\{U_{ik}U_{jk}\}_{k\in\N}$ is an i.i.d. sequence of random
  variables, for all $i$, $j\in\{1,\dots,N\}$ and, due to
    the Strong Law of Large Numbers \cite[p.~6]{AWV00}, when $N\to \infty$ we have
  \begin{align*}
    \frac{1}{N}\sum_{k=1}^{N} U_{ik}U_{jk} \xrightarrow{\text{a.s.}}
    \E[U_{i1}U_{j1}] =
    \begin{cases}
      \sigma^{2}, & \text{if } i=j, \\
      0, & \text{if } i\ne j,
    \end{cases}
  \end{align*}
  where $\E[\cdot]$ denotes the expected value operator. Then,
  \begin{align}\label{eq:UUT}
    \frac{1}{N}U U^{\transpose} \xrightarrow{\text{a.s.}} \sigma^{2}I
    \ \ \text{ as }  N\to \infty.
  \end{align} 
  Next, consider the function $f:\R^{mT\times mT}\to \R^{mT\times n}$, $
  Y \mapsto YC_T^{\transpose} ( C_T Y C_T^{\transpose})^{\dag}$.
  Note that $f(Y)$ is continuous at $Y=\alpha I$, $\alpha>0$,\footnote{In fact, since $\text{Rank}(C_T Y C_T^\transpose)=\text{Rank}(C_TC_T^\transpose)$ for any positive definite $Y$, it holds $\lim_{k\to \infty} (C_T Y_k C_T^\transpose)^\dag = (\alpha\, C_T C_T^\transpose)^\dag$ for any sequence of positive definite matrices $\{Y_k\}_{k\in\N}$ such that $\lim_{k\to \infty} Y_k = \alpha I$~\cite[p.~238]{AB-TNEG:03}.} and $f(\alpha I)=C_T^\transpose (C_T C_T^\transpose)^\dag = C_T^\dag$ \cite[p.~49]{AB-TNEG:03}.
  To conclude, we employ the Continuous Mapping Theorem  \cite[Theorem 2.3(iii)]{AWV00} and \eqref{eq:UUT} to obtain, as $N\to\infty$,
  \begin{align*}
    U X^\dag &= U(C_T U)^{\dag}= \frac{1}{N}U U^{\transpose}C_T^{\transpose} \left(C_T \frac{1}{N}U U^{\transpose} C_T^{\transpose}\right)^{\dag}\\
             & = f \left(\frac{1}{N}UU^{\transpose} \right) \xrightarrow{\text{a.s.}} f \left(\sigma^{2} I\right) = C_T^{\dag} . 
  \end{align*}
\end{IEEEproof}
Theorem \ref{thm: equivalent expression 2} contains a data-driven
expression of the minimum-energy control input for a linear system,
which does not rely on the estimation of the system matrices or the
controllability matrix. As we show in the next section, the expression
\eqref{eq: solution M input} is not only conceptually simpler than the
classic Gramian-based expression of the minimum-energy control input
and our other data-driven expressions \eqref{eq: U min 1} and
\eqref{eq: alternative expression}, but it is also numerically more
reliable as it requires a smaller number of operations. Yet,
differently from \eqref{eq: U min 1} and \eqref{eq: alternative
  expression}, the expression \eqref{eq: solution M input} coincides
with the minimum-energy control only asymptotically in the number of
experiments, and assuming that the entries of the input matrix
  $U$ are zero-mean i.i.d. random variables with nonzero
  finite~variance.

\begin{remark}{\bf \emph{(Data-driven minimum-energy control inputs
      when $x_0 \neq 0$)}}\label{remark: nonzero initial state}
  When $x_0 \neq 0$, the computation of the minimum-energy control
  input to reach $\subscr{x}{f}$ is more involved, as the unknown
  matrix $A$ and vector $x_0$ enter the relation \eqref{eq: data
    relation}.\footnote{Notice that the term $A^{T}x_0$ remains
    unknown even if the exact value of $x_0\ne 0$ is known. Thus
    knowledge of $x_0$ does not modify the expressions we obtain when
    $x_0\ne 0$ is treated as an unknown~variable.} Yet, under a mild
  assumption on the experimental inputs $U$, minimum-energy inputs can
  still be computed with a finite number of experiments. To see this,
  consider the problem
  \begin{align}\label{eq: min energy problem initial state}
    \begin{array}{ll}
      \min\limits_{\alpha} &  \| U \alpha \|_2^2 ,\\[1em]
      \,\text{s.t.} 
                           & \subscr{x}{f} = X \alpha \  \text{ and } \  1 = \vectorones^\transpose \alpha,
    \end{array}
  \end{align}
  Assume that the matrix $U$ is full row rank, and that there exists a
  vector $w$ such that $Uw = 0$ and $\vectorones^\transpose w \neq
  0$. The first assumption guarantees that there exists $\alpha^{*}$
  such that $u^{*}=U\alpha^{*}$, and thus the computation of the
  minimum-energy control for any final state $\subscr{x}{f}$
  (cf. Theorem 2.1). The second assumption ensures that there exists
  $\alpha^{*}$ satisfying $1=\1^{\transpose}\alpha^{*}$, which allows
  us to correctly reconstruct the term $A^{T}x_0$ from
  $X$.\footnote{These assumptions can always be satisfied by
    properly designing the experimental inputs, or by running
    sufficiently many random experiments.
  } In fact, let
  $ \alpha^* = U^\dag u^* + w\,{(1 - \vectorones^\transpose U^\dag
    u^*)}/{(\vectorones^\transpose w)}, $ and notice that
  $u^* = U\alpha^*$, where $u^*$ is the minimum-energy control input
  to reach $\subscr{x}{f}$. Further, using~\eqref{eq: data relation}
  and $1 = \vectorones^\transpose \alpha^*$, we have 
  $X \alpha^* = \sum_{i = 1}^N X_i \alpha_i^* = A^T x_0 \sum_{i =
    1}^N \alpha_i^* + C_T \sum_{i = 1}^N \alpha_i^* U_i = A^T x_0 +
  C_T u^* = \subscr{x}{f} .  $ Then, similarly to the proof of Theorem
  \ref{thm: feasibility}, a solution to \eqref{eq: min energy problem
    initial state} determines the minimum-energy input.

  To solve the minimization problem \eqref{eq: min energy problem
    initial state}, let
  $\bar X = [X^\transpose \; \vectorones]^\transpose$ and
  $\subscr{\bar x}{f} = [\subscr{x}{f}^\transpose \;
  1]^\transpose$. Then, similarly to  Theorem \ref{thm: feasibility},
  we obtain
  $
    \alpha^* = \bar X^\dag \subscr{x}{f} - K (UK)^\dag U \bar X^\dag
    \subscr{\bar x}{f},
  $
  where $K = \Ker(\bar X)$, and
  \begin{align}\label{eq: input non zero initial conditions}
    u^* = (I - UK (UK)^\dag) U \bar X^\dag \subscr{\bar x}{f} .
  \end{align}
  Because the matrix $U$ is required to have a nontrivial
  null-space, a sufficient number of linearly-independent non-optimal
  experiments for the computation of the minimum-energy control input
  to any arbitrary final state is $mT + 1$.
  
  Finally, from the above reasoning and the proof of Theorem \ref{thm:
    equivalent expression} and Theorem \ref{thm: equivalent expression
    2}, the minimum-energy input \eqref{eq: input non zero initial
    conditions} can be written equivalently as
  $ u^* = (\bar X U^\dag)^\dag \subscr{\bar x}{f} = U \bar X^\dag
  \subscr{\bar x}{f}$,
  where the last equality holds asymptotically for any choice of
  inputs satisfying the assumptions in Theorem \ref{thm:
    equivalent expression 2}.~\oprocend
\end{remark}

\begin{figure}
 \hspace{0.15cm}
\includegraphics[scale=1]{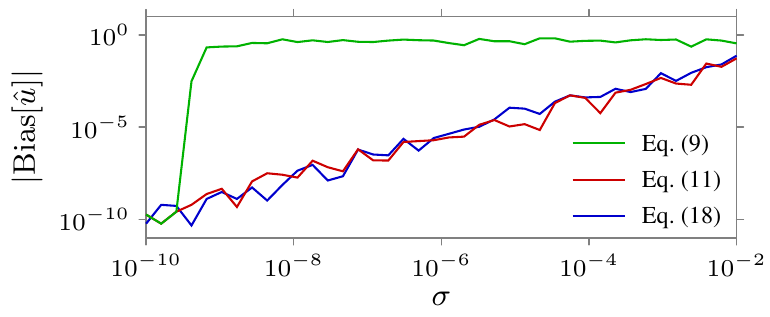}
\caption{This figure shows the magnitude of the bias of the data-driven expressions \eqref{eq: U min 1}, \eqref{eq: alternative expression}, and \eqref{eq: solution M input} as a function of the standard deviation of the noise $\sigma$. We choose $A={\left[\protect\begin{smallmatrix}-0.8 & 0 & 0 \protect\\ 2 & 0.1 & 0 \protect\\ 0.2 & 1 & 0.5\protect\end{smallmatrix}\right]}$, $B={\left[\protect\begin{smallmatrix}1\protect\\ 0 \protect\\ 0\protect\end{smallmatrix}\right]}$, $\subscr{x}{f}={\left[\protect\begin{smallmatrix}0.3\protect\\ 1 \protect\\ 0.5\protect\end{smallmatrix}\right]}$, $T=8$, and $N=10$. The entries of $U$ and $X$ have been chosen randomly and then corrupted by i.i.d.~Gaussian noise with zero mean and standard deviation $\sigma$. The bias has been computed as the average over 100 noise~realizations.}
\label{Fig:noise}
\end{figure}

\begin{remark}{\bf \emph{(Data-driven expressions with noisy
      data)}}\label{remark: noise}
  Let the measurements of the input $u_i$ and the final state $x_i$ be
  corrupted by noise. Let
  $\tilde{U}=[u_{1}+w_{1}\ \cdots\ u_{N}+w_{N}]$ and
  $\tilde{X}=[x_{1}+v_{1}\ \cdots \ x_{N}+v_{N}]$ be the matrices
  obtained by concatenating all noisy measurements. The data-driven
  estimates \eqref{eq: U min 1}, \eqref{eq: alternative expression},
  and \eqref{eq: solution M input} computed from the noisy data
  $(\tilde U, \tilde X)$ are typically biased. To see this, consider
  the system $x(t+1)=ax(t)+u(t)$, $a\in \real$, $x_0=0$, and
  $T=N=1$. In this simple case, expressions \eqref{eq: U min 1},
  \eqref{eq: alternative expression}, and \eqref{eq: solution M input}
  are equivalent and, assuming that $x_{1}+v_{1}\ne 0$, read as
  $\hat u = \frac{u_{1}+w_{1}}{x_{1}+v_{1}}x_{\text{f}}$. If $w_{1}$
  and $v_{1}$ are independent random variables uniformly distributed
  in $[-\varepsilon, \varepsilon]$, with $0<\varepsilon<|u_{1}|$, it
  holds
  \begin{align*}
    \mathrm{Bias}[\hat u] &=\mathbb{E}_{w_{1},v_{1}}[\hat u]-u^{*} =
    \mathbb{E}_{v_{1}}\left[\frac{u_{1}}{u_{1}+v_{1}}\right]x_{\text{f}}-x_{\text{f}}
    \\
    &=
    \left[\frac{1}{2\varepsilon}u_{1}\ln\left(\frac{u_{1}+\varepsilon}{u_{1}-\varepsilon}\right)-1\right]x_{\text{f}},
  \end{align*}
  where $\mathbb{E}_{z}[\cdot]$ denotes the expected value with
  respect to $z$. It can be shown that, if $u_{1}$ and $x_{\text{f}}$
  are nonzero, the previous equation vanishes only in the limit
  $\varepsilon\to0$. This implies that all data-driven expressions in
  this simple case are biased. When $n>1$, a quantitative
  characterization of the bias (and covariance) of the data-driven
  expressions appears to be difficult, due to the presence of
  pseudoinverse operations. However, numerical simulations with
  i.i.d.~normally distributed noise (see also Fig.~\ref{Fig:noise})
  suggest that (i) all data-driven expressions are biased in the case
  of noisy measurements, (ii) the magnitude of the bias is
  proportional to the standard deviation $\sigma$ of the noise for
  \eqref{eq: alternative expression} and \eqref{eq: solution M input},
  while  it increases rapidly as $\sigma$ grows and sets to a
  constant value for \eqref{eq: U min 1}.
  \oprocend
 \end{remark}
 
  \begin{remark}{\bf \emph{(Data-driven expressions with output measurements)}}\label{remark: output}
  Consider the system
\begin{align*}
  x(t+1) = A x(t) + Bu(t),\ \ y(t) = Cx(t),
\end{align*}
where $C\in\real^{p\times n}$, and assume that for each experimental
input $u_{i}$, $i\in \until{N}$, we can measure the output of the
system at time $T$, namely, $y_{i}=Cx_{i}$. Let
$Y=[y_{1}\ \cdots\ y_{N}]\in\real^{p\times N}$ be the matrix
concatenating all output measurements, and assume that the system is
output controllable in $T$ steps. That is, the $T$-steps output
controllability matrix $C_{O,T}=[CB \ \ CAB \ \ \cdots\ \ CA^{T-1}B]$
has full row rank \cite{EK-PS:64}. The minimum-energy input to reach
the output $y_{\text{f}}\in\real^{p}$ in $T$ steps is
$u^{*}=C_{O,T}^{\dagger}(y_{\text{f}}-CA^{T}x_0)$. All results
discussed in this paper apply to the case of output control after
substituting $X$ and $x_{\text{f}}$ with $Y$ and $y_{\text{f}}$,
respectively.  \oprocend
 \end{remark}

\begin{figure*}[t]
\begin{center}
\includegraphics[scale=1]{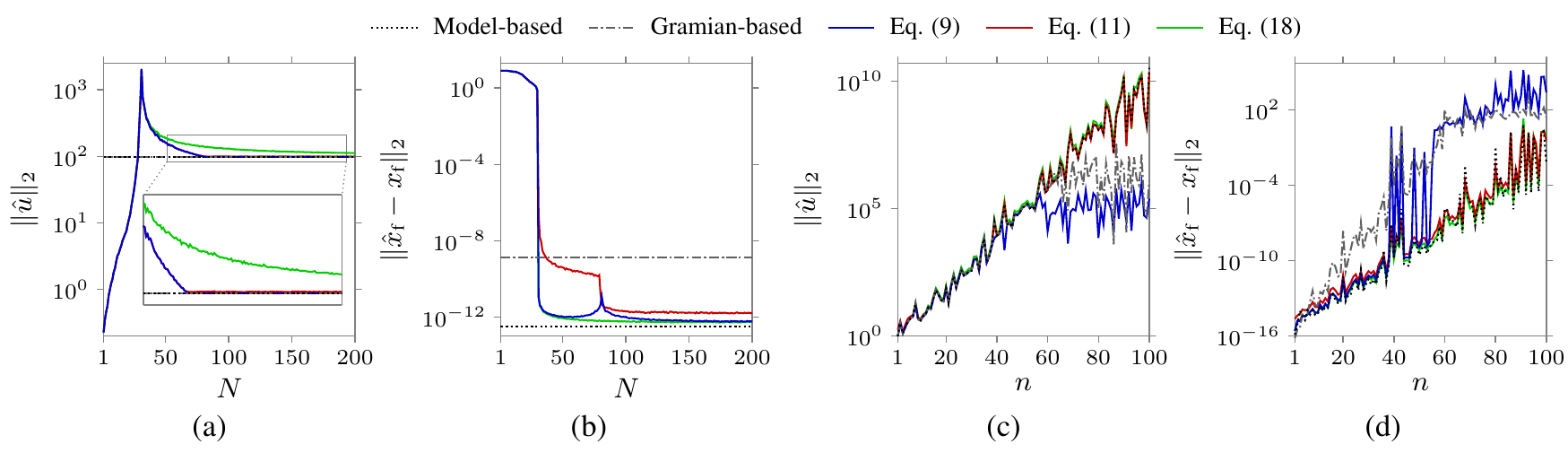}
\caption{Fig.~\ref{Fig:numerical_analysis_full}(a)-(b) show the norm of the control input
  $\hat{u}$ computed via the model-based formula \eqref{eq:
    model-based} (dotted line), via inversion of the controllability
  Gramian (dashed line), and via the data-driven expressions
  \eqref{eq: U min 1}, \eqref{eq: alternative expression}, and
  \eqref{eq: solution M input} extended for $x_0\ne 0$ as in Remark
  \ref{remark: nonzero initial state} (colored lines), and the
  corresponding error in the final state, as the number of data $N$
  varies. We choose $n=20$, $m=2$, and $T=40$. The matrix $A$
  (possibly unstable) has been populated with random i.i.d. normal
  entries and then normalized by $\sqrt{n}$, the entries of $B$,
  $x_0$ and $\subscr{x}{f}$ have been chosen randomly according to a
  normal distribution. The curves represent the average over 100
  experiments with data pairs $(x_i,u_i)$, where $u_i$ has random
  i.i.d. normal entries, and
  $x_i=C_Tu_i$. Fig.~\ref{Fig:numerical_analysis_full}(c)-(d) show the
  norm of the inputs $\hat{u}$ computed as above, and the
  corresponding errors in the final state, as a function of the system
  dimension $n$. We choose $m=2$, $T=n$, and $N=mT+20$. The matrices
  $A$ and $B$ have been generated as above. The curves represent the
  average over 1000 experiments with data $(x_i,u_i)$ and states
  $x_0$, $\subscr{x}{f}$ generated as above. All the computations
  have been carried out using standard built-in \texttt{Matlab 2016b} linear
  algebra routines.}
\label{Fig:numerical_analysis_full}
\end{center}
\end{figure*}

\section{Numerical analysis}

What remains unclear from the previous analysis is the benefit, if
any, in collecting a large number of control experiments. We next show
that increasing the number of control experiments can improve the
numerical reliability and accuracy of computing minimum-energy control
inputs.  

In Fig.~\ref{Fig:numerical_analysis_full} we compare the numerical
performance of the model-based expressions of the minimum-energy
controls $u^* = C_T^\dag \subscr{x}{f}$ and
$u^* = C_T^{\transpose} W_T^{\dag} \subscr{x}{f}$ (Gramian-based),
with our data-driven expressions in \eqref{eq: U min 1}, \eqref{eq:
  alternative expression}, and \eqref{eq: solution M input}. In
particular, in Fig.~\ref{Fig:numerical_analysis_full}(a)-(b) we plot
the norm of the control inputs and the numerical errors in reaching
the final state $\subscr{x}{f}$, for all strategies and as a function
of the number $N$ of control experiments. Here, we focus on a
  ``worst-case'' analysis and choose a small input dimension ($m=2$),
  since a large value of $m$ certainly improves the conditioning of
  all expressions. Fig.~\ref{Fig:numerical_analysis_full}(a) shows
that the norm of the data-driven control inputs \eqref{eq: U min 1}
and \eqref{eq: alternative expression} equals its minimum value when
$N\ge mT$ (as predicted by Theorems \ref{thm: feasibility} and
\ref{thm: equivalent expression}), whereas the norm of the data-driven
input \eqref{eq: solution M input} converges to its minimum value only
asymptotically (as predicted by Theorem \ref{thm: equivalent
  expression 2}). Fig.~\ref{Fig:numerical_analysis_full}(b) shows
that, for sufficiently large $N$, the final state reached by the three
data-driven control strategies is almost as close to $\subscr{x}{f}$
as the one computed via the model-based formula
$u^* = C_T^\dag \subscr{x}{f}$, and considerably closer to
$\subscr{x}{f}$ than the state reached by the Gramian-based control
input, with expressions \eqref{eq: U min 1} and \eqref{eq: solution M
  input} being the most accurate, showing that the computation of the
minimum-energy control input via our data-driven expression is as
reliable as the computation of the input based on the exact knowledge
of the system matrices, and numerically more reliable than the
model-based Gramian formula. Instead, in
Fig.~\ref{Fig:numerical_analysis_full}(c)-(d) we plot the norm of the
control inputs obtained through the different strategies described
above and their corresponding errors in the final state as a function
of the system dimension $n$. As expected, the accuracy of the
Gramian-based control input deteriorates rapidly as $n$
increases. Yet, surprisingly, the data-driven expressions of the
minimum-energy control inputs remain accurate for systems of
considerably larger dimension. Further, the data-driven control
\eqref{eq: solution M input} yields the smallest error in the final
state among the three data-driven strategies. This could be due to the
simpler form of \eqref{eq: solution M input}, which requires the
computation of only one pseudoinverse, or to the fact that the energy
of \eqref{eq: solution M input} reaches the minimum value only
asymptotically in $N$. Finally,
Fig. \ref{Fig:numerical_analysis_full}(c)-(d) show that expression
\eqref{eq: U min 1} becomes numerically unreliable for smaller values
of the system dimension compared to \eqref{eq: alternative expression}
and \eqref{eq: solution M input}. This is likely because of the
additional computations~in~\eqref{eq: U min 1}.

\section{Conclusion and future work}\label{sec: conclusion}
In this paper we derive data-driven expressions of open-loop
minimum-energy control inputs for linear systems. Leveraging linearity
of the dynamics, we show that such optimal controls can be learned
from a finite number of control experiments, without knowing or
reconstructing the system matrices, and where the control experiments
are conducted with non-optimal and arbitrary inputs. We derive
  three different data-driven expressions of minimum-energy controls:
  while \eqref{eq: alternative expression} appears to be the simplest
  exact data-driven expression, \eqref{eq: U min 1} constitutes a
  radically different and new way of computing minimum-energy
  controls, and highlights several geometric connections between the
  minimum-energy solutions and the experimental data, and \eqref{eq:
    solution M input} provides a simple way of computing a family of
  data-driven, sub-optimal, minimum-energy controls.  We further
illustrate that our data-driven expressions of the minimum-energy
inputs are simpler and numerically more reliable than the classic
Gramian-based expression, 
especially when the dimension of the system increases.

The results of this paper support the intriguing idea of combining
model-based control methods with data-driven techniques, showing that
this new framework has the potential to considerably increase the
reliability and effectiveness of the two parts alone. This paper also
creates several directions of future research, including the
  extension to closed-loop, noisy, and model predictive control
  problems.

\bibliographystyle{unsrt}
\bibliography{BIB}

\end{document}